\documentclass{amsart}
\usepackage{amssymb} 
\newtheorem{thm}{Theorem}[section]

\newtheorem{cor}[thm]{Corollary}
\newtheorem{lem}[thm]{Lemma}
\newtheorem{pro}[thm]{Proposition}
\newtheorem{DEF}[thm]{Definition}
\newtheorem{exa}[thm]{Example}
\newtheorem{qe}[thm]{Question}
\newtheorem{rem}[thm]{Remark}
\numberwithin{equation}{section}

\begin{document}

\title[]{ configuration of nilpotent groups and isomorphism}%
\author{A. Abdollahi, A. Rejali and A. Yousofzadeh}%

\email{a.abdollahi@math.ui.ac.ir}%
\email{rejali@sci.ui.ac.ir}%
\email{a.yousofzade@sci.ui.ac.ir}%
\thanks{ }
\date{}%
\begin{abstract}
The concept of configuration was first introduced by Rosenblatt
and Willis to give a condition for amenability of groups.  We show
that if $G_1$ and $G_2$ have the same configuration sets and $H_1$
is a normal subgroup of $G_1$ with abelian quotient, then there is
a normal subgroup $H_2$ of $G_2$ such that
$\frac{G_1}{H_1}\cong\frac{G_2}{H_2}.$ Also configuration of
FC-groups and isomorphism is studied.
\end{abstract}
\maketitle
\date{}%
\maketitle

\section{Introduction} The notion of configuration for a finitely
generated group $G$ was introduced  in \cite{RW7}. It was shown
that amenability of $G$ is characterized by its configuration
equations. In \cite{ARW1}, the authors investigated some
properties of groups which can be characterized by configurations
and studied the question of whether $G$ is determined up to
isomorphism by its configurations.

 The configurations of $G$ are defined in terms of
finite generating sets and finite partitions of $G$. A
configuration corresponding to a generating sequence $\frak g
=(g_1,\dots,g_n)$  and a partition $\mathcal E=\{E_1,\dots,E_m \}$
of $G$ is an $(n+1)-$tuple $C=(c_0,\dots,c_n)$, where $1\leq
c_i\leq m$ for each $i$, such that there is $x$ in $G$ with $x\in
E_{c_0}$ and $g_ix\in E_{c_i}$ for each $1\leq i\leq n$. The set
of all configurations corresponding to the pair $(\mathfrak g,
\mathcal E)$ will be denoted by $Con(\mathfrak g, \mathcal E)$.
$$Con(G)=\{Con(\mathfrak g, \mathcal E)|\ \frak g\ a\  finite\ generating\
sequence\
of\
G,\ \mathcal E\ a\ finite\ partition \}  $$ The set of all
configuration sets of $G$ is denoted by $Con(G).$ We write
$G_1\approx G_2$ and $G_1$ is called configuration equivalent with
$G_2$, if $Con(G_1)=Con(G_2)$. Clearly  the relation $\approx$ is
an equivalence relation on the class of finitely generated groups.

The configuration $C=(c_0,\dots,c_n)$ may be described
equivalently as a labeled tree. The tree has one vertex of degree
$n$, labeled by $c_0$. Emanating  from this vertex are edges
labeled $1,\dots,n$, and the other vertex of the $i$-th edge is
labeled $c_i$.  When the generators are distinct, this tree is a
subgraph of the Cayley graph of the finitely generated group
$G=\langle g_1,\dots,g_n\rangle$. The edges labels indicate which
generator gives rise to the edge and the vertices labels show
which set of the partition $\mathcal E$ the vertex belong to. From
this perspective the configuration set $Con(\mathfrak g, \mathcal
E)$ is a set of rooted trees having height 1. This finite set
carries information about $G$ and the present paper addresses some
properties of $G$ which can be recovered from such information.
We prove that if the group $G$ has the same configuration set with
the direct product of groups $\mathbb Z^n\times F,$ where $F$ is a
finite group, then $G$ is isomorphic with $\mathbb Z^n\times F$.
This generalizes the result of Abdollahi, Rejali and Willis
\cite{ARW1} saying that  if $G_1$ is an abelian group with the
same configuration sets as $G_2$, then they are isomorphic.

In \cite{ARW1} the authors tried to generalize the latter result
for nilpotent groups. It was shown that if $G_1$ is a finitely
generated nilpotent group of class $c$ having the same
configuration sets with $G_2$, then $G_2$ is a nilpotent group of
class $c$. In this paper we show that if $G_1$ is a torsion free
nilpotent group of Hirsch length $h$, then so is $G_2$.

We do not know whether the configuration equivalent and
isomorphism relations are the same for the class of nilpotent
groups. We will  prove it  for certain nilpotent groups.

Let $G$ be a finitely generated nilpotent group. Then $G$ is
isomorphic with a subgroup of direct product of $T_r(n, \mathbb
Z)$, and a finite nilpotent group, for some positive integer $n$,
where $T_r(n, \mathbb Z)$ denotes the group of all upper
triangular matrices with integer entries and all diagonal entries
equal 1. It is well-known that  $T_r(n, \mathbb Z)$ is a finitely
generated torsion free nilpotent group of class $n-1$ and Hirsch
length $\frac{n(n-1)}{2}$ (see \cite{R6}).
  We study finitely generated torsion free nilpotent
groups $G$ of class 2 such that $\frac{G}{T}\cong \mathbb Z^n$ and
$T\cong \mathbb Z^m$, where $T=\tau(G)$ is the isolator of the
commutator subgroup of $G$, and this class of groups will be
denoted by $\mathcal I(n,m)$ (see \cite{S8}). It is shown that if
$G\in \mathcal I(n,m)$ has the same configuration set with $H$,
then $H \in \mathcal I(n,m)$. Therefore
$\frac{G}{\tau(G)}\cong\frac{H}{\tau(H)}$ and
$\tau(G)\cong\tau(H).$ For $n=3$ and $m=2$, groups $G$ and $H$
will be isomorphic.

It is shown in  \cite{ARW1} that if a finitely generated group
$G_1$ has the same configuration set with $G_2$, and $G_1$ has a
normal subgroup $H_1$ with finite index, then $G_2$ has a normal
subgroup $H_2$ so that $\frac{G_1}{H_1}\cong\frac{G_2}{H_2}$. For
a group $G$, we denote by $\mathcal{F}(G)$ the set of all
isomorphism types of finite quotients of $G$. Thus the latter
result says that the set of isomorphism types of finite quotients
of two configuration equivalent groups $G_1$ and $G_2$ coincide
i.e., $\mathcal F(G_1)=\mathcal F(G_2)$. We show that one can
extend finite index property to abelian quotient property. In fact
$\mathcal A(G_1)=\mathcal A(G_2),$ where $\mathcal A(G)$ is the
set of all isomorphism types of
abelian quotients of $G$.\\

Let  $G$ be a finitely generated infinite simple group.  Then
$\mathcal F({\mathbb Z}_2)=\mathcal F(G \times \mathbb{Z}_2)$,
however $\mathbb{Z}_2$ and $G\times \mathbb{Z}_2$ are not
configuration equivalent (see \cite{ARW1}). This example shows
that the relation of having the same configuration is strictly
stronger than the relation of having the same finite quotients.
  Several authors characterized groups with the same finite
quotient sets, under certain conditions (see e.g., [8] and [9]).
In \cite{S8} it was shown that
 two polycyclic-by-finite  groups $G_1$ and $G_2$ have the same
finite quotient sets if and only if $\widehat{G_1}\cong
\widehat{G_2}$, where $\widehat{G}$ is the profinite completion of
$G$.   Also in \cite{war13}   cancellation of groups is studied
 and under certain conditions it is shown that $\mathcal F(G_1)=\mathcal
 F(G_2)$ if and only if $\mathbb Z\times G_1\cong \mathbb Z\times
 G_2$. By using this
 result one can show that if $G_1$ and $G_2$
 are finitely generated FC-groups with the same configuration sets,
 then $G_1$ is isomorphic with a subgroup of $G_2$ of finite
 index, where by  an FC-group we mean a group  in which  each
 element has finitely many conjugates. In particular if $G_1$ is
 a finitely generated nilpotent  FC-group, then so is $G_2$.
 We do not know whether it is true without the nilpotence condition.

 It is to be noted that a finitely generated group $G$ is an
 FC-group if and only if $\frac{G}{Z(G)}$ is finite, where
 $Z(G)$ is the center of $G$ (see \cite{T10}). We prove that
if $G_1$ and $G_2$ are two finitely generated nilpotent FC-groups
having the same configuration sets, then
$\frac{G_1}{Z(G_1)}\cong\frac{G_2}{Z(G_2)}$ and $Z(G_1)\cong
Z(G_2)$. So configuration equivalent nilpotent FC-group of
nilpotency class $2$ have equivalent upper central series. We do
not know if this result is true for general nilpotent groups. If
$G$ is a finitely generated   FC-group, then $G$ is isomorphic
with a subgroup of $\mathbb Z^n\times F$, for some finite group
$F$, (see \cite{T10}). We show that configuration equivalence and
isomorphism are the same for $\mathbb Z^n\times F$.

Let $T_r(n, \mathbb Z)=\langle g_{i,j}\ :\ 1\leq i <j\leq
n\rangle,$ in which $g_{i,j}=I+E_{i,j}$, where $E_{i,j}$ is the
matrix with entry 1 in the $(i,j)-$th position and other entries
are zero. Also $I$ is the $n\times n$ identity matrix. One can
show that each matrix $A=(a_{i,j})\in T_r(n, \mathbb Z)$ has a
unique representation as
$$
A=(g_{1,n}^{a_1,n}g_{2,n}^{a_2,n}\cdots
g_{n-1,n}^{a_{n-1},n})(g_{1,n-1}^ {a_1,n-1}g_{2,n-1}^{a_2,n-1}
\cdots g_{n-2,n-1}^{a_{n-2},n-1})\cdots g_{1,2}^{a_{1,2}}.$$
 Furthermore $g_{i,j}g_{k,t}=g_{k,t}g_{i,j}$ and
$g_{k,t}g_{i,k}=g_{i,k}g_{i,t}^{-1}g_{k,t}$ for all $i\neq t$ and
$j\neq k$.

\section{Configuration equivalence of nilpotent groups  and
isomorphism }

In this section we are looking for group theoretical properties
$\mathcal{P}$ satisfying the condition of the question posed in
\cite{ARW1}: Which group properties are translated by
configuration equivalence? For a group $G$, we denote by
$\mathcal{F}(G)$ and $\mathcal{A}(G)$ the set of isomorphism types
of finite quotients and abelian quotients of $G$, respectively. We
first mention a result in this direction which was proved in
\cite{ARW1}.
\begin{thm}{\rm(See Lemma 6.3 of \cite{ARW1}
)}\label{M-F} Let $G_1$ and $G_2$ be  finitely generated groups
such that $G_1 \thickapprox G_2$. Then
$\mathcal{F}(G_1)=\mathcal{F}(G_2)$.
\end{thm}
We actually have further for two configuration equivalent  groups,
indeed the sets of their abelian quotients are the same. We need
the following lemma to prove the latter.
\begin{lem}\label{F-A}
Let $G$ be a finitely generated abelian group and $L$ be a finite
abelian group. If for every prime $p$, there exists a subgroup
$N_p$ of $G$ such that $\frac{G}{N_p}\cong \mathbb{Z}_p^n\times
L$, then $G$ has a subgroup $N$ such that $\frac{G}{N}\cong
\mathbb{Z}^n\times L$.
\end{lem}
\begin{proof} Suppose that
\begin{center}
$G\cong \mathbb{Z}^m\times
\mathbb{Z}_{{p_1}^{\alpha_1}}\times\ldots\times
\mathbb{Z}_{{p_k}^{\alpha_k}}$,
\end{center}
where $p_1,\dots,p_k$ are prime numbers. Then every finite factor
of $G$ is isomorphic to
\begin{center}
$ \mathbb{Z}_{t_1}\times\ldots\times\mathbb{Z}_{t_m}\times
\mathbb{Z}_{{p_1}^{\beta_1}}\times\ldots\times
\mathbb{Z}_{{p_k}^{\beta_k}}$,
\end{center}
for some nonnegative integers
$t_1,\dots,t_m,\beta_1,\ldots,\beta_k$, where
$\beta_i\leq\alpha_i$. For, if $M$ is a subgroup of $G$, then
there exist a generating set $\{g_1,\dots,g_m,h_1,\dots,h_k\}$ of
$G$ and positive integers $d_1,\dots,d_{m+k}$ such that $\langle
g_i\rangle\cong\mathbb{Z}, \langle
h_i\rangle\cong\mathbb{Z}_{{p_i}^{\alpha}}$,
\begin{center}
$G=\langle g_i\rangle\times\cdots\times\langle
g_m\rangle\times\langle h_1\rangle \times\cdots\times\langle
h_k\rangle $
\end{center}
and $M=\langle
g_1^{d_1},\dots,g_m^{d_m},h_1^{d_{m+1}},\dots,h_k^{d_{m+k}}\rangle
$. Therefore
\begin{center}
$\frac{G}{M}\cong\frac{\langle g_1\rangle
\times\cdots\times\langle g_m\rangle \times\langle h_1\rangle
\times\cdots\times\langle h_k\rangle } {\langle g_1^{d_1}\rangle
\times\cdots\times\langle g_m^{d_m}\rangle \times\langle
h_1^{d_{m+1}}\rangle \times\cdots\times\langle
h_k^{d_{m+k}}\rangle }$.
\end{center}
This gives the result. Now let $p$ be a prime number greater than
$|L|$ and $p_1^{\alpha_1}\cdots p_k^{\alpha_k}$. Then by
assumption and the first part of the proof, there exist
nonnegative integers $t_1,\dots,t_m$ and $\beta_1,\ldots,\beta_k$
with $\beta_i\leq\alpha_i$ such that
\begin{center}
$\mathbb{Z}_{t_1}\times\cdots\times\mathbb{Z}_{t_m}\times\mathbb{Z}_{p_1^{\beta_1}}\ldots
\times\mathbb{Z}_{p_k^{\beta_k}}\cong\mathbb{Z}_{p}^n\times L$.
\end{center}
Now by the uniqueness of the elementary factors of a finitely
generated abelian group, we must have $n\leq m$,
\begin{center}
$\mathbb{Z}_{l_1}\times\ldots\times\mathbb{Z}_{l_n}\cong\mathbb{Z}_p^n$
and
$\mathbb{Z}_{r_1}\times\ldots\times\mathbb{Z}_{r_{m-n}}\times\mathbb{Z}_{p_1^{\beta_1}}\times\ldots
\times\mathbb{Z}_{p_k^{\beta_k}}\cong L$,
\end{center}
where $\{t_1,\ldots,t_m\}=\{l_1,\ldots,l_n,r_1,\ldots,r_{m-n}\}$.
Therefore
\begin{center}
$\frac{\mathbb{Z}^m\times\mathbb{Z}_{p_1^{\alpha_1}}\times\ldots\times
\mathbb{Z}_{p_k^{\alpha_k}}} {\{0\}^n\times
r_1\mathbb{Z}\times\ldots\times
r_{m-n}\mathbb{Z}\times\mathbb{Z}_{p_1^{\alpha_1-\beta_1}}\times\ldots\times
\mathbb{Z}_{p_k^{\alpha_k-\beta_k}}}\cong\mathbb{Z}^n\times L$.
\end{center}
This completes the proof.
\end{proof}
\begin{pro}
Let $G$ and $H$ be finitely generated groups such that
$G\thickapprox H$. Then $\mathcal{A}(G)=\mathcal{A}(H)$.
\end{pro}
\begin{proof} By Lemma \ref{F-A}, it is enough to prove that if
$\frac{G}{N} $ is an (infinite) abelian factor of $G$, then
$\frac{G}{N}\cong \frac{H}{M}$ for some normal subgroup $M$ of
$H$. Since $G$ is finitely generated,
$\frac{G}{N}\cong\mathbb{Z}^n\times L$, where $n$ is a positive
integer and $L$ is a finite abelian group. Since
$\mathbb{Z}^n\times L$ has finite factors $\mathbb{Z}_p^n\times L$
for any prime number $p$, by Theorem 3.1, there exists a normal
subgroup $M_p$ of $H$ such that
$\frac{H}{M_p}\cong\mathbb{Z}_p^n\times L$. Now let
$B:=\bigcap_{p\  prime}M_p.$ Then $\frac{H}{B}$ is an abelian
group having $\mathbb{Z}_p^n\times L$ for all prime number $p$ as
a factor group. It follows from  Lemma \ref{F-A} that
$\frac{H}{B}$ has a factor group isomorphic to $\mathbb{Z}^n\times
L$, as required. This completes the proof.
\end{proof}

\begin{pro}\label{pol-tor}
Let $G$ and $H$  be two  polycyclic groups such that $G\approx H$.
If $G$ admits a normal series with infinite cyclic factors, then
$H$ is torsion free.
\end{pro}
\begin{proof} By Theorem \ref{M-F} $\mathcal{F}(G)= \mathcal{F}(H)$.  Now
by using  Lemma 1.5 of \cite{T11} (see also chapter 10 of
\cite{S8}),
 $H$ is
torsion free.
\end{proof}
In a polycyclic group $G$ the number of infinite factors in a
cyclic series, which is known as the Hirsch length, is independent
of the series and is an invariant of $G$. Denote by $h(G)$ the
Hirsch length of a polycyclic group $G$.

\begin{pro}\label{lem-nil-tor}
Let $G$ and $H$ be two finitely generated groups such that
$G\approx H$. If $G$ is  torsion free nilpotent of class $c$ with
the Hirsch length $h$, then so is $H$.
\end{pro}
\begin{proof}
By  Corollary 5.2 of \cite{ARW1}, $H$ is also nilpotent of class
$c$. Since every finitely generated torsion-free nilpotent group
admits a normal (central) series with infinite cyclic series,
Proposition \ref{pol-tor} implies that $H$ is torsion-free. Now
one can prove by induction on $c$, that the Hirsch lengths of $G$
and $H$ coincide.
\end{proof}
 Let $G$ be a finitely generated torsion-free nilpotent
group of class 2. We denote by   $T=\tau(G)$ the isolator of $G'$
in $G$, i.e.  the set of all elements $x\in G$ such that $x^s\in
G'$ for some non zero integer $s$. The set $T$ is a central
subgroup of $G$ and $T/G'$ is the torsion subgroup of $G/G'$. Then
there exist positive integers $m$ and $n$ such that
$\frac{G}{T}\cong \mathbb{Z}^n$ and $T\cong \mathbb{Z}^m$. The
class of all such groups $G$ shall be denoted by
$\mathcal{I}(n,m)$ (see \cite{S8}, {P. 260}).
\begin{lem}\label{I(m,n)}
Let $G\in \mathcal{I}(n,m)$ and $G\approx H$ for some finitely
generated  group $H$.  Then $H\in \mathcal{I}(n,m)$.
\end{lem}
\begin{proof} By Proposition \ref{lem-nil-tor}, $H$ is a finitely
generated torsion-free nilpotent group of class 2. Suppose that
$H\in\mathcal{I}(n',m')$.
 We have
$\frac{G}{\tau(G)}\cong \mathbb{Z}^n$ and $\frac{H}{\tau(H)}\cong
\mathbb{Z}^{n'}$. By Proposition 3.3, there exists a normal
subgroup $K$ of $H$ such that $\frac{H}{K}\cong \mathbb{Z}^n$. If
$x\in \tau(H)$, then $x^s \in H'$ for some non zero integer $s$.
Since $H/K$ is abelian, $H'\leq K$ and so $x^s \in K$. As $H/K$ is
torsion-free, we have that $x\in K$. Thus $\tau(H)\leq K$. It
follows that $n'\geq n$. Similarly $n\geq n'$ and so $n=n'$. \\
By Proposition \ref{lem-nil-tor}, $h(G)=h(H)$. On the other hand,
$h(G)=h(\tau(G))+h(G/\tau(G))$ and $h(H)=h(\tau(H))+h(H/\tau(H))$.
Therefore $h(\tau(G))=h(\tau(H))$.
 It follows that $m=m'$, as
required.
\end{proof}
\begin{thm}\label{IIq}
Let $G\in \mathcal{I}(n,m)$ and $G\approx H$ for some finitely
generated  group $H$. If $(n,m)\in\{(3,2),(3,3),(n,1)\}$, then
$G\cong H$.
\end{thm}
\begin{proof}  By Lemma \ref{I(m,n)}, $H\in \mathcal{I}(n,m)$. It follows from
Theorem \ref{M-F}, $\mathcal{F}(G)=\mathcal{F}(H)$. Now  results
stated in p. 265 of \cite{S8} complete the proof.
 \end{proof}
 \begin{cor}
Let $G$ be a finitely generated group such that
$Con(G)=Con(Tr(3,\mathbb{Z}))$. Then $G\cong Tr(3,\mathbb{Z})$.
 \end{cor}
 \begin{proof}
It is easy to see that $Tr(3,\mathbb{Z})\in \mathcal{I}(2,1)$. Now
Theorem \ref{IIq} completes the proof.
 \end{proof}
The following question is natural:

\begin{qe} \em
Let $G$ be a finitely generated group. Does
$Con(G)=Con(Tr(n,\mathbb{Z}))$ imply $G\cong Tr(n,\mathbb{Z})?$
\end{qe}

\begin{rem} \em In \cite{RY} the authors gave a positive answer to the above
question with another type of configuration, say two-sided
configuration.
\end{rem}
\section{Configuration equivalent for FC-groups }

It is interesting to know if being FC-group is conserved by
equivalence of configuration. We answer this question in  some
special cases.

It is shown that for each finitely generated FC-group $G$,
 the torsion subgroup $Tor(G),$ and  commutator subgroup $G'$ of
$G$ are finite and the factor group $\frac{G}{Tor(G)}$ is an
abelian torsion free group \cite{T10}. It has been shown that if
two groups are configuration equivalent, then the isomorphism
classes of their finite quotients are the same \cite{ARW1}. For
groups with finite commutator there is an other interesting
result:

\begin{thm}\label{funda}
Let $G$ and $H$ are two finitely generated groups with finite
commutator subgroups. Then the following conditions are
equivalent:

(a) $\mathcal{F}(G)=\mathcal{F}(H), $

(b) $G\times \mathbb{Z}\cong H\times\mathbb{Z}.$
\end{thm}
\begin{proof}
See [10, Theorem 2.1].
\end{proof}
\begin{DEF}
For any group $G$, the $non\ cancellation\ set,\ \chi(G)$ is the
set of isomorphism classes of groups $H$ such that
$G\times\mathbb{Z}\cong H\times \mathbb{Z}$.
\end{DEF}
It is easy to see that for a finite or abelian group $G$,
$\chi(G)$ is trivial.

The following theorem is proved in \cite{wi14}.
\begin{thm}\label{fund}
Let $G$ and $H$ be any groups with finite commutator subgroups. If
$G$ is infinite and $\chi(G)$ is trivial than $\chi(G\times H)$ is
trivial.
\end{thm}

So for example if $G$ is an infinite abelian group and  $H$ is
arbitrary, then $\chi(G\times H)$ is trivial.

\begin{thm}
Let $Con(G)=Con(\mathbb{Z}^n\times F),$ where $F$ is a finite
group and $G$ is an FC-group. Then $G\cong\mathbb{Z}^n\times F$.
\end{thm}
\begin{proof}
By Theorem (\ref{funda}) $\mathbb{Z}\times G\cong \mathbb{Z}\times
\mathbb{Z}^n\times F.$ But $\chi(\mathbb{Z}^n\times F)$ is trivial
by Theorem (\ref{fund}). Therefore $G\cong\mathbb{Z}^n\times F$.
\end{proof}

By  a direct proof we show in the next theorem  that the condition
of being FC-group for $G$ is superfluous. This states the main
result of this section.

 For a non empty set $A$, we denote by $\chi_B$ the characteristic
 function of a subset $B$ of $A$, which is defined as
 $$\chi_B(a)=\begin{cases} 1 & \text{if} \; a\in B\\
 0 & \text{if} \; a\in A\backslash B. \end{cases}$$

\begin{thm}
Let $Con(G)=Con(\mathbb{Z}^n\times F),$ where $F$ is a finite
group. Then $G\cong\mathbb{Z}^n\times F$.
\end{thm}
\begin{proof} Suppose that $F=\{x_1,\ldots ,x_l\}$ where $x_1=e$.
Let $\Sigma=\{\sigma_1,\ldots,\sigma_{3^n}\}$ be the set of all
functions from $\{1,\ldots,n\}$ to $\{0,1,-1\}$. Put
$E_{\sigma,j}=\sigma(1)\mathbb{N}\times\ldots\times
\sigma(n)\mathbb{N}\times\{x_j\},$ for any $\sigma\in \Sigma$.
Then $\mathcal{E}=\{E_{\sigma,j}: \ \sigma\in\Sigma ,\
j\in\{1,\ldots,l\}\}$ is a partition of $\mathbb{Z}^n\times F.$
Let   $g=(g_1,\ldots ,g_{l+n})$ where $g_1=(1,0,\ldots
,0,e)$,\ldots ,$g_n=(0,0,\ldots ,1,e),$ for $1\leq i\leq n$ and
$g_{n+j}=(0,0,\ldots ,0,x_j)$ for $1\leq j\leq l$. Then
$\mathbb{Z}^n\times F=\langle g_1,\ldots ,g_{l+n}\rangle$. Let
also $\pi:\{1,\ldots ,l\}\times\{1,\ldots ,l\}\rightarrow
\{1,\ldots ,l\}$ be defined by $\pi (i,j)=k$ whenever  $x_i x_j
=x_k$.

We refine the  partition $\mathcal{E}$ to a new partition
$\mathcal{E}'$ such that

\noindent if $\sigma= \chi_{\{i\}}$ , then
$E_{\sigma,1}=A_i\bigcup B_i$ where $A_i=\{g_i\}$ and
$B_i=E_{\sigma,1}\backslash A_i$.

\noindent if $\sigma=\chi_{\{i,j\}}$ then
$E_{\sigma,1}=A_{i,j}\bigcup B_{i,j}$, in which
$A_{i,j}=\{g_ig_j\}$ and $B_{i,j}=E_{\sigma,1}\backslash A_{i,j}$.

 Let $(f,\mathcal{K})$ be the configuration pair of $G$, in which
$\mathcal{K}=\{K_{\sigma,j}: \ \sigma\in\Sigma ,\
j\in\{1,\ldots,l\}, \ (\sigma,j)\neq (\chi_{\{i\}},1) ,\ \
(\sigma,j)\neq
(\chi_{\{i,k\}},1)\}\bigcup\{C_i,D_i,C_{i,k},D_{i,k}\ :\ 1\leq
i,k\leq n\}$ is a partition and $f=(f_1,\ldots ,f_{l+n})$ is a
generating set of $G$,
  such that
$Con(f,\mathcal{K})=Con(g,\mathcal{E}).$

After writing the members of $Con(f,\mathcal{K}),$ we have

\noindent(I)  $f_iK_{\sigma,j}\subseteq K_{\sigma'_i,j}$ where
$\sigma(i)=0$ and $\sigma'_i(i)=1$  and $\sigma(t)=\sigma'_i(t)$
for all $t\neq i$,

\noindent(II) $f_iK_{\sigma,j}\subseteq K_{\sigma,j}$ where
$\sigma(i)=1,$

\noindent(III)  $f_iK_{\sigma,j}\subseteq K_{\sigma,j}\bigcup
K_{\sigma''_i,j}$ where $\sigma''_i(i)=0$ and $\sigma(i)=-1$ and
$\sigma(t)=\sigma''_i(t)$ for all $t\neq i$,

\noindent(IV) $f_{n+i}K_{\sigma,j}= K_{\sigma,\pi(i,j)}$, for
$1\leq i \leq l$ and $1\leq j \leq l.$

 Now let $\sigma(i)=0$ then
for each $j\in\{1,\ldots ,l\}$ and each $m\in \mathbb{N}$ we have
$f_i^m K_{\sigma,j}\subseteq K_{\sigma'_i,j}$ and so
$\bigcup_{m=1}^{\infty} f_i^m K_{\sigma,j}\subseteq
K_{\sigma'_i,j}$. On the other hand $f_i( K_{\sigma,j} \bigcup
K_{\sigma'_i,j})= K_{\sigma'_i,j}$.

Consider the case $\sigma(i)=-1$. For each $j\in\{1,\ldots ,l\}$,
$K_{\sigma,j }=f_i^{-1}(K_{\sigma,j}\bigcup K_{\sigma''_i,j})$. It
is routine  to compute that
$\bigcup_{m=1}^{\infty}f_i^{-m}K_{\sigma''_i,j}\subseteq
K_{\sigma,j}$

Finally if $K_1=K_{\sigma_{0},1}$ where $\sigma_0\equiv 0$ and
$1\in K_1$  and $\sigma\in\Sigma$, then
$f_if_jK_1=f_iC_j=C_{i,j}=f_jC_i=f_jf_iK_1$ and,
\begin{center}
\begin{equation}
\bigcup_{m_1=1}^{\infty}\cdots \bigcup_{m_n=1}^{\infty}
f_1^{m_1\sigma(1)}\cdots f_n^{m_n\sigma(n)}f_{n+j}K_1 \subseteq
K_{\sigma,j}. \label{asl}
\end{equation}
\end{center}
 We claim that  each element of $G$ belongs to
$f_1^{m_1}\cdots f_n^{m_n}f_{n+j}K_1$, for some integers
$m_1,\ldots ,m_n$ and some $j\in\{1,\ldots ,l\}$.

Suppose that $K_1\neq\{1\}$ and set $L_{\sigma_0,1}=\{1\}$ and
$L_2=K_{\sigma_0,1}\backslash\{1\}$. Then if $L_a=K_a$, for each
$a$,
\begin{center}
$\mathcal{L}=\{L_a\ :\ a\neq \sigma_{0,1}\}\bigcup L_1\bigcup L_2$
\end{center}

 is a partition of $G$. So by
hypotheses there exist a configuration pair $(\mathcal{T},h)$ of
$\mathbb{Z}^n\times F$ such that
$Con(\mathcal{T},h)=Con(\mathcal{L},f)$. Let $h_s=(a_s^1,\ldots
,a_s^n,x_{j_s})$ for some $a_s^i\in\mathbb{Z}$ and $x_{j_s}\in F$,
where $1\leq s \leq n+l$. We claim that $a_{n+i}^t=0$ for each
$1\leq i \leq l$ and $1\leq t \leq n$. If they are  not, since the
vectors $(a_1^1,\ldots ,a_1^n),\ldots ,(a_n^1,\ldots ,a_n^n)$ and
$(a_{n+i}^1,\ldots ,a_{n+i}^n)$ are $\mathbb{Q}$ linearly
dependent in $\mathbb{Q}^n$ and so in $\mathbb{Z}^n$, there exist
some integers $k_1,\ldots ,k_n$ and $k_i$ not together all zero
such that $h_1^{k_1}\cdots h_n^{k_n}h_{n+i}^{k_i}=(0,\ldots
,0,x_l)$ for some $x_l\in F$ and so $h_1^{k_1}\cdots
h_n^{k_n}h_{n+i}^{k_i}$ must have a finite order. But it is
possible only if $k_1=k_2=\cdots =k_n=0$, by (\ref{asl}).
Therefore
\begin{center}
$h_{n+i}=(0,\ldots ,0,x_{j_{n+i}}).$
\end{center}

Since $h_{n+1},\ldots,h_{n+l}$  are $l$ distinct members,
($h_{n+i}\in K_{\sigma_0,i}$) we have $\{x_{j_{n+1}},\ldots ,
x_{j_{n+l}}\}=F.$

Now consider that
\begin{center}
$\mathbb{Z}^n\leq \langle (a_1^1,\ldots ,a_1^n),\ldots
,(a_n^1,\ldots ,a_n^n) \rangle.$
\end{center}

Consequently each element of $\mathbb{Z}^n\times F $ has the form
$h_1^{k_1}\cdots h_n^{k_n}h_{n+j}.$ Also $1\in T_1$ and
$h_1^{k_1}\cdots h_n^{k_n}h_{n+j}\in T_{\sigma,j}$, in which
$\sigma(i)=sgn(k_i),$ is defined $1,-1$ and $0$ if $k_i$ is
positive, negative and zero respectively. So
\begin{center}
$\mathbb{Z}^n\times F\subseteq (\bigcup _{T\in
\mathcal{T}}T)\backslash T_2.$
\end{center}
 Thus $T_2=\emptyset.$

This contradiction causes that $K_1=\{1\}$ i.e. each element of
$G$ is of the form $f_1^{k_1}\cdots f_n^{k_n}f_{n+j},$ and

\begin{center}
$K_{\sigma,j}=\{ f_1^{k_1\sigma(1)}\cdots
f_n^{k_n\sigma(n)}f_{n+j} \ , \ \ k_1,\ldots ,k_n\in\mathbb{N}\}$
\end{center}

Note that the expression of elements are unique. For if
 \begin{center}
$x=f_1^{k_1}\cdots f_n^{k_n}f_{n+j}=f_1^{s_1}\cdots
f_n^{s_n}f_{n+t},$

\end{center}

then there exists a $\sigma\in \sum$ such that $x\in
K_{\sigma,j}\bigcap K_{\sigma,t}$, and so $t=j$. Therefore

\begin{center}
$f_1^{k_1}\cdots f_n^{k_n}=f_1^{s_1}\cdots f_n^{s_n},$

$f_1^{k_1-s_1}\cdots f_n^{k_n}=f_2^{s_2}\cdots f_n^{s_n}.$
\end{center}
So by (\ref{asl}), $\sigma(1)=0$, which implies $k_1=s_1$. By
induction $k_i=s_i$ for every $i\in\{1,\ldots ,n\}.$
 Define
\begin{center}
$\phi:G\rightarrow \mathbb{Z}^n\times F$,

$ \phi(f_1^{m_1}\cdots f_n^{m_n}f_{n+j})=(m_1,\ldots ,m_k,x_j),$
\end{center}
Evidently it is a one to one  function onto $\mathbb{Z}^n\times
F$.

It is also a group isomorphism since
\begin{center}
$\phi(f_1^{k_1}\cdots f_n^{k_n}f_{n+j}f_1^{k'_1}\cdots
f_n^{k'_n}f_{n+t})= \phi(f_1^{k_1}\cdots f_n^{k_n}f_1^{k'_1}\cdots
f_n^{k'_n}f_{n+\pi(j,t)})=\phi(f_1^{k_1+k'_1}\cdots
f_n^{k_n+k'_n}f_{n+\pi(j,t)}) =(k_1+k'_1,\ldots
,k_n+k'_n,x_{\pi(i,j)})$.
\end{center}
This completes the proof.

\end{proof}

 In the following  $n=n(G)=n_1n_2n_3$, in which $n_1$ is the exponent of
the torsion subgroup of $G,$ $n_2$ is the exponent of
$Aut(Tor(G))$ and $n_3$ is the exponent of torsion group of
$Z(G).$ (The exponent of a group is the least $n\in \mathbb N$
such that $x^n=1$ for all $x\in G$ or $\infty$ if such $n$ does
not exist.)

Let $G$ have finite commutator subgroup and $n=n(G)$. If $G\approx
H$, then by [11, Theorem 4.2], $H$ is isomorphic to a subgroup $L$
of $G$ of finite index in $G$ that $[G:L]$ is relatively prime to
$n$. Now let $[G:L]\equiv \pm 1$ $mod(n)$. Then $H\cong G$  by
[11, Theorem 4.3].

\begin{thm}
Let $G$ be a finitely generated nilpotent FC-group and $G\approx
H$. Then $H$ is a nilpotent FC-group too. Furthermore $Z(G)\cong
Z(H)$ and $\frac{G}{Z(G)}\cong \frac{H}{Z(H)}$. Also $G$ and $H$
have the same Hircsh length.
\end{thm}
\begin{proof}
Since $G$ is an FC-group, $\frac{G}{Z(G)}$ is finite and
$|\frac{G}{Z(G)}|=n$ say. So $xy^n=y^nx$ for each $x,y$ in $G$.
Therefore by \cite{ARW1} this semigroup law is satisfied by the
group $H.$ Thus $\frac{H}{Z(H)}$ is a periodic finitely generated
group. It is also nilpotent. On the other it is known that every
periodic finitely generated nilpotent group is finite. Therefore
$\frac{H}{Z(H)}$ will be finite and so $H$ is an FC-group by
\cite{T10}. Now  using Theorem \ref{funda}, $G\times
\mathbb{Z}\cong H\times\mathbb{Z}.$ Therefore $Z(G)\times \mathbb
Z\cong Z(H)\times \mathbb{Z}.$ But $Z(G)$ and $Z(H)$ are abelian
and so $Z(G)\cong Z(H),$ by [9, Theorem 7]. Also we can easily see
that $\frac{G}{Z(G)}\cong\frac{H}{Z(H)}$. On the other hand since
$\frac{G}{Z(G)}$ and $\frac{H}{Z(H)}$ are finite and $Z(G)\cong
Z(H),$   by  [6, p.16] $G$ and $H$ have the same Hircsh lengths.
\end{proof}
 The following example shows that configuration equivalence  and
 isomorphism is the same  for a group which is neither an FC-group,
  nor a nilpotent one.

\begin{exa}\label{Dinfty}
Let $G$ be a finitely generated group such that $G\thickapprox
D_{\infty}$, where $D_{\infty}=\langle x,y\;|\; x^2=y^2=1\rangle$
is the infinite dihedral group. Then $G\cong D_{\infty}$.
\end{exa}
\begin{proof}
  Consider the following configuration pair
$(\mathfrak{g},\mathcal{E})$ of $G$: $\mathfrak{g}=(x,y)$ and
$\mathcal{E}=\{E_1,\dots,E_5\}$, where
$$E_1=\{1\},\ E_2=\{x\},\  E_3=\{y\},$$
$$E_4=\{g_1g_2\cdots g_n\;|\;\ g_i\in\{x,y\}, n\in \mathbb{N},n>1
,\ i\in\{1,\dots,n\}, g_1=x,g_{i}\neq g_{i+1}\}, \;\text{and}$$
$$E_5=\{g_1g_2\cdots g_n\;|\; \ g_i\in\{x,y\}, n\in \mathbb{N},n>1
,\ i\in\{1,\dots,n\}, g_1=y, g_{i}\neq g_{i+1}\}.$$
 Evidently $\mathcal{E}$ is a partition of $D_\infty$. We have
$Con(\mathfrak{g},\mathcal{E})=\{c_1,\dots,c_7$\}, where $$
\begin{matrix}
   c_1=(1,2,3),  & c_2=(2,1,5),&c_3=(3,4,1),& c_4=(4,5,5), \\
    c_5=(4,3,5), & c_6=(5,4,2),& c_7=(5,4,4).  \\
  \end{matrix}$$
Since  $Con(G)=Con(D_{\infty})$, there exists an ordered
generating set $\mathfrak{h}=(h_1,h_2)$ and a partition
$\mathcal{F}=\{F_1,\dots,F_5\}$ of $G$ such that
$Con(\mathfrak{g},\mathcal{E})=Con(\mathfrak{h},\mathcal{F})$ and
$1\in F_1$. Thus
  $$\begin{matrix}
    h_1F_1=F_2,& h_1F_2=F_1,& h_2F_1=F_3,& h_2F_3=F_1,\\
h_1F_4=F_3\cup F_5,& h_2F_5=F_2\cup F_4, &h_2(F_2\cup F_4)= F_5, &
h_1(F_3\cup F_5)=F_4.&
  \end{matrix}$$
 Therefore $h_1^2F_1=h_2^2F_1=F_1$ and for $n>1$ and $t_1,\dots,t_n
 \in \{h_1,h_2\}$ with $t_i\neq t_{i+1}$ we have $n\geq 1,$

$$t_1\cdots t_nF_1\subseteq \begin{cases} F_4 &\text{if}\;\; t_1=h_1
\\ F_5 &\text{if}\;\; t_1=h_2 \end{cases},\ \ \ \ \ \ (t_i\in \{h_1,h_2\}).$$
Now let
 $$ \begin{matrix}
L_1=\{1\},& L_2=F_1\backslash\{1\},& L_3=F_2,\\
 L_4=F_3,& L_5=F_4,& L_6=F_5.
     \end{matrix}$$
 and $\mathcal{L}=\{L_1,\ldots,L_6\}$. If
$L_2\neq\varnothing$, then $\mathcal{L}$ is a partition of $G$.
Therefore we have $Con(\mathfrak{h},\mathcal{L})=
\{d_1,\ldots,d_{10}\}$, where
  $$\begin{matrix}
d_1=(1,3,4),& d_2=(2,3,4),& d_3=(3,1,6),& d_4=(3,2,6),& d_5
=(4,5,1),\\ d_6 =(4,5,2),& d_7=(5,4,6),& d_8=(5,6,6),& d_9 =(6,5,3
),&  d_{10}=(6,5,5).
  \end{matrix}$$
Thus there exist an ordered generating set
$\mathfrak{g}'=(g_1,g_2)$ and a  partition
$\mathcal{K}=\{K_1,\ldots,K_6\}$ of $D_{\infty}$ such that
$Con(\mathfrak{g}',\mathcal{K})=\{d_1,\ldots,d_{10} \}$ and $1\in
K_1$. Thus

\begin{eqnarray}\label{www}
g_1(K_1\cup K_2)=K_3,& g_1K_3=K_1\cup K_2,\nonumber\\ g_1(K_6\cup
K_4)=K_5,&
g_1K_5=K_4\cup K_6,\\
g_2(K_1\cup K_2)= K_4,& g_2K_4= K_1\cup K_2, \nonumber\\
g_2(K_3\cup K_5)= K_6,& g_2K_6= K_3\cup K_5.\nonumber
\end{eqnarray}

Note that an element of $D_{\infty}$ has finite order if and only
if its length as product of $a$ and $b$ is odd.

 Three different cases are possible:

Case I) $g_1$ (respectively $g_2$) has even order and $g_2$
(respectively $g_1$) has odd order. In this case $g_1g_2$ must
have an odd length and consequently finite order, which is
impossible since it is easy to see that $(g_1g_2)^n\in K_5\cup
K_6$ for each non zero integer $n.$

Case II) Both $g_1$ and $g_2$ have even lengths. In this case the
group $\langle g_1,g_2\rangle$ can not have any non identity
element of finite order and it is a contradiction too, since
$D_{\infty}=\langle g_1,g_2\rangle.$

\smallskip

Case III) Both $g_1$ and $g_2$ have odd lengths. Thus they must be
of order $2$. Also by (\ref{www}) each element of $D_{\infty}$ has
a unique form of  products of $g_1$ and $g_2.$ On the other hand $
1\in K_1, g_1\in K_3$ and $g_2\in K_4$. We have also for $n\geq
1,$

$$t_1\cdots t_nK_1\subseteq \begin{cases} K_5 &\text{if}\;\; t_1=g_1
\\ K_6 &\text{if}\;\; t_1=g_2 \end{cases},\ \ \ \ \ \ (t_i\in \{g_1,g_2\}).$$

In other word, $D_{\infty}\subseteq K_1\cup K_3\cup K_4\cup
K_5\cup K_6.$ So $K_2$ should be empty.
  This implies that  $F_1=\{1\}$ and so
$h_1^2=h_2^2=1$. Since $h_1$ and $h_2$ are not in $F_1$,
$o(h_1)=o(h_2)=2$. By Proposition 4.3 of \cite{ARW1}, $G$ is
infinite and so $G$ is an infinite group generated by two
involutions. Hence $G\cong D_{\infty}$, as required.
\end{proof}
 In \cite{ARW1} the authors proved that if a finitely generated group $G$
satisfies a semigroup law and $G\approx H$, then so does $H$. In
the previous  example $D_{\infty}$ is a finitely presented group
and we show their result for generators in this case. Now the
following question is natural.

\begin{qe}
Let $G\approx H$ and the generating sequence $\mathfrak
g=(g_1,\ldots ,g_n)$ of $G$ satisfies the semigroup law $w(\frak
g)=1.$ Does there exist a generating sequence $\mathfrak h$ of $H$
satisfying $w(\mathfrak h)=1?$
\end{qe}
\noindent{\bf{Acknowledgements.}} The first two authors's research
was partially supported by the Center of Excellence for
Mathematics at University of Isfahan. The first author's research
was in part supported by a grant from IPM (No. 87200118).

\vspace{.5cm}

\noindent{
\address{Ali Rejali, Department of Mathematics, University of Isfahan,
Isfahan 81744, Iran}}

\vspace{.5cm}

\noindent{
 \address{Alireza Abdollahi, Department of Mathematics, University
 of Isfahan, Isfahan 81744-73441, Iran and School of Mathematics, Institute
for Research in Fundamental Sciences (IPM), P.O.Box: 19395-5746,
Tehran, Iran. }}

\vspace{.5cm}

\noindent{\address{Akram Yousofzadeh, Department of Mathematics,
University of Isfahan, Isfahan 81744, Iran}}

\bigskip

\end{document}